# Digits of pi: limits to the seeming randomness II

Paula Nataniela Roba, Karlis Podnieks

University of Latvia
19 Raina Blvd., Riga, LV-1586, Latvia
paulanatanielar@gmail.com, karlis.podnieks@lu.lv

**Abstract.** According to a popular belief, the decimal digits of mathematical constants such as π behave like statistically independent random variables, each taking the values 0, 1, 2, 3, 4, 5, 6, 7, 8, and 9 with equal probability of 1/10. If this is the case, then, in particular, the decimal representations of these constants should tend to satisfy the *central limit theorem* (CLT) and the *law of the iterated logarithm* (LIL).

The paper presents the results of a direct statistical analysis of the decimal representations of 12 mathematical constants with respect to the central limit theorem (CLT) and the law of the iterated logarithm (LIL). The analyzed constants include π, e, the Euler–Mascheroni constant (γ), Apéry's constant (ζ(3)), quadratic irrationals such as the square roots of 2 and 3, and the golden ratio. Additionally, natural logarithms – *log 2, log 3*, and *log 10*, Catalan's constant (G), and the lemniscate constant ($\bar{\omega}$) were examined.

The first billion digits of each constant were analyzed, with ten billion digits examined in the case of π. Within these limits, no evidence was found to suggest that the digits of these constants satisfy CLT or LIL.

## 1. Introduction

The decimal digits of $\pi, e, \sqrt{2}$ and other "naturally occurring" mathematical constants are widely believed to behave like statistically independent random variables, each taking the values 0, 1, 2, 3, 4, 5, 6, 7, 8, and 9 with equal probability of 1/10.

Where does the idea come from that the digits of these constants "must" behave randomly? Could this belief stem from a vague intuition that the way "naturally occurring" mathematical constants are defined is fundamentally "orthogonal" to the principles governing the decimal (or any other base) expansions of real numbers?

However, what does this popular belief actually mean in precise terms? Could it have any "practical" consequences? Does it even make sense? After all, random variables take random values, whereas, for example, the decimal digits of π form a single, deterministic, and computable sequence $s_{10}(\pi) = 3.14159\ldots$!

We propose setting aside concerns about the "randomness" of these deterministic sequences of digits and instead compare their behavior with the probabilistic model of "statistically independent random variables that take the values 0–9 with equal probability 1/10", as considered in probability theory.

Namely, to generalize from decimal digits to an arbitrary base $q$, let us consider an infinite sequence of statistically independent random variables $X_1^{(q)}, \ldots X_n^{(q)}, \ldots$ taking values from the set $Q_q = \{0, 1, \ldots, q-1\}$ with equal probability 1/q. The corresponding probability space $Q_q^N$ consists of all infinite sequences of base-$q$ digits, where each sequence can be interpreted as a real number. The natural probability measure on this space is defined by assigning a probability of $q^{-n}$ to each cylinder set of sequences (a cylinder set consists of all sequences that extend a given initial string of digits $d_1, \ldots, d_n$).

Of course, $s_{10}(\pi) \in Q_{10}^N$, and this fact does not depend on the probability measure defined on $Q_q^N$.

However, a fundamental issue with such infinite models is that nontrivial hypotheses about them cannot be verified empirically. Instead, theoretical analysis is required – confirming or refuting hypotheses by proving theorems (for example, within Zermelo–Fraenkel set theory, ZFC). In the worst case, some hypotheses may turn out to be undecidable in ZFC, much like the famous Continuum Hypothesis.

For example, if we were to successfully define the "randomness property" of sequences as some subset $R_q \subset Q_q^N$, how could we determine whether $s_{10}(\pi) \in R_{10}$ or not? Most likely, this could not be decided empirically; the only way to resolve it would be to prove either $s_{10}(\pi) \in R_{10}$, or $s_{10}(\pi) \notin R_{10}$ as a theorem.

In the absence of a precise definition of "randomness," many relaxed versions have been proposed, known as "randomness tests." These tests formulate properties that must be satisfied by all infinite "random" sequences.

For instance, "random" sequences must not be periodic. Since $\pi$ has been proven to be irrational, it follows as a theorem that, for any base $q$, the base-$q$ digits of $\pi$ are not periodic. Unfortunately, beyond this, little has been proven about the *infinite* sequence of digits of $\pi$.

In particular, we are still unable to prove that all digits 0, 1 ,…, 9 appear in the decimal expansion of $\pi$ with equal frequency 1/10 – a hypothesis proposed by Émile Borel in (1909). All we can do is attempt to test this hypothesis empirically, examining a finite number of digits.

Peter Trueb (2016) conducted such an empirical analysis using his record-breaking computation of $n \approx 22.4592 \cdot 10^{12}$ decimal digits of $\pi$. For such $n$, in the probabilistic model, the expected variance of digit frequencies is $\sigma^2 = 0.1(1-0.1)/n \approx 4.01 \cdot 10^{-15}$, corresponding to a standard deviation of $\sigma = 0.63 \cdot 10^{-7}$.

As shown in Figure 1 of his paper, for the first $22.4592 \cdot 10^{12}$ digits of $\pi$, the observed frequencies deviate from the expected 1/10 by no more than approximately $-1.5\sigma$ to $+1.9\sigma$. The reported observed variance of the frequencies is about $5.0 \cdot 10^{-15}$—approximately 25% higher than the variance expected under the probabilistic model. Should this be regarded as a confirmation of Borel's hypothesis?

## 2. Related work

Thanks, ChatGPT, for the assistance.

To the best of our knowledge, the most recent advances are presented in three papers. On page 261 of his 1909 paper, Borel generalized the definition of "normality" from single digits to any "groupements" of digits. He proposed taking any string $J$ of $k$ digits and sliding it along the sequence of the first $n$ digits of the number $C$, counting matches in $count_q(C,n,J)$ and thus obtaining the frequency:

$$freq_q(C,n,J) = \frac{1}{n} count_q(C,n,J).$$

If, for each $n$ and any $k$-digit string $J$,

$$\lim_{n \to \infty} freq_q(C,n,J) = \frac{1}{q^k},$$

then Borel calls the number $C$ "normal in base $q$."

This "sliding pattern" definition of normality may seem unusual because pattern windows are allowed to overlap. However, in a "truly random" sequence, no sliding "groupement" is expected to deviate significantly from the binomial distribution with match probability $p_k = q^{-k}$.

**Trueb (2016)** reports the results of normality tests using his record-breaking dataset of 22.4 trillion decimal (and 18.6 trillion hexadecimal) digits of $\pi$. He has tested Borel's normality for $q=10$ and $q=16$, and for k=2 and k=3:

"The frequencies of all sequences up to length 3 in the first 22'459'157'718'361 decimal and 18'651'926'753'033 hexadecimal digits of $\pi$ are found to be consistent with the hypothesis of $\pi$ being a normal number in base 10 and base 16."

**Sevcik (2018)** reports the application of fractal analysis:

"Fractal analysis is used here to show that $\pi$'s digit sequence corresponds to a uniformly distributed random succession of independent decimal digits, and that these properties get clearer as the number of digits in the series grows towards infinity; $10^9$ digits were tested in this work. This is consistent with the hypothesis that $\pi$ is normal."

"Infinity is far away. The findings in this study are based on the first $10^9$ initial digits of $\pi$; this limit was set by the computer power available to the author [Sevcik]. ... the most interesting aspect of this study is that the more series digits you consider, the stronger its conclusions get. … Fractal analysis presented here indicates that $\pi$ is normal in base 10."

**Plouffe (2022)** reports extensive pattern searching tests:

"A series of large-scale tests were performed on the first 1000 billion digits of the number $\pi$. … The purpose of these tests is to detect possible patterns. Other tests were made on ten mathematical constants to 1 billion decimal places."

"The first 1000 billion decimal places of π were used and 10 constants with 1 billion decimal places: $\sqrt{2}, \sqrt{3}, \phi, e, \ln(2), \ln(10), \gamma$, Lemniscate constant, $\zeta(3)$ and Catalan."

"The result of the tests is negative, no pattern was found."

"The current record for calculation [of π] is 10^14 decimal places and to do a test you need quite a bit more than a single pc or two and a few hard drives. In my [Plouffe's] opinion, only Google can do a test at the same level. … There is no example of a natural constant in base 2 or 10 or even in any base that exhibits any pattern."

## 3. Central Limit Theorem (CLT)

The advanced results reported in Section 2, consider, in a sense, somewhat derived settings.

We will follow the classical straightforward approach. In the probabilistic model described above, we consider an infinite sequence of statistically independent random variables $X_1^{(q)}, .. X_n^{(q)}, ...$, each taking values from the set $Q_q = \{0, 1, ..., q-1\}$ with equal probability *1/q*.

We will analyze the sums of the first *n* digits, namely:

$$S_n^{(q)} = \sum_{i=1}^{n} X_i^{(q)}.$$

The expected value $\mu$ and variance $\sigma^2$ of each digit $X_i^{(q)}$ can be computed as follows:

$$\mu = \sum_{j=0}^{q-1} \frac{1}{q} j = \frac{q-1}{2} ; \sigma^2 = \sum_{j=0}^{q-1} \frac{1}{q}(j-\mu)^2 = \frac{q^2-1}{12},$$

and its third absolute moment (for even *q*) is given by:

$$\alpha_3 = \sum_{j=0}^{q-1} \frac{1}{q} |j-\mu|^3 = \frac{q(q^2-2)}{32}.$$

For decimal digits (*q=10*) this yields:

$$\mu = 4.5 ; \sigma^2 = 8.25 ; \sigma \approx 2.872 ; \alpha_3 = 30.625.$$

Hence, the expected value of the sum $E(S_n^{(q)}) = \mu n$ , its variance $Var(S_n^{(q)}) = \sigma^2 n$, its standard deviation is $\sigma \sqrt{n}$ , and its third absolute momentum $M_3(S_n^{(q)}) = \alpha_3 n$.

Aiming for the Central Limit Theorem (CLT), let us consider the corresponding centered and normalized deviation of the sum $S_n^{(q)}$ from its expected value $\mu n$:

$$d_n^{(q)} = \frac{S_n^{(q)} - \mu n}{\sigma \sqrt{n}}.$$

The variables $d_n^{(q)}$ take their values in the same probabilistic space as the variables $X_1^{(q)}, \ldots, X_n^{(q)}, \ldots$, namely, in $Q_q^N$. Let us denote by $F_n^{(q)}(x)$ the cumulative distribution function (CDF) of of $d_n^{(q)}$, i.e., the probability:

$$F_n^{(q)} = P(d_n^{(q)} \leq x).$$

According to CLT, as $n \to \infty$, the function $F_n^{(q)}(x)$ converges to the cumulative distribution function of the standard normal distribution, namely, the function:

$$\Phi(x) = P(X \leq x) = \frac{1}{\sqrt{2\pi}} \int_{-\infty}^{x} e^{-\frac{t^2}{2}} dt \ .$$

The convergence rate of $F_n^{(q)}(x)$ to $\Phi(x)$ can be estimated by applying the Berry-Esseen inequality: for all $n$ and $x$,

$$\left| F_n^{(q)}(x) - \Phi(x) \right| \leq 0.4748 \frac{\alpha_3}{\sigma^3 \sqrt{n}}.$$

This version of the inequality, with the constant 0.4748,was proved by Shevtsova (2011).

For decimal digits ($q=10$) this would mean that for all $n$ and $x$:

$$\left| F_n^{(10)}(x) - \Phi(x) \right| \leq 0.4748 \frac{30.625}{2.872^3 \sqrt{n}} < \frac{0.61}{\sqrt{n}}.$$

Thus, for example, starting with $n = 10^6$, the probability $P(d_n^{(10)} \leq x)$ will differ from $\Phi(x)$ less than $10^{-3}$, and we can safely assume that $P(d_n^{(10)} \leq x) \approx \Phi(x)$.

Now, let us turn to the digits of mathematical constants.

Consider base-$q$ representations (q>1) of some mathematical constant $C$.

Let us denote by $D_q(C, i)$ the i-th digit in the base-$q$ representation of the number C, and by $S_q(C, n)$ – the sum of the first $n$ digits in this representation, namely:

$$S_q(C, n) = \sum_{i=1}^{n} D_q(C, i) \ .$$

And, as in the probabilistic model, let us consider the corresponding centered and normalized deviation:

$$d_q(C, n) = \frac{S_q(C, n) - \mu n}{\sigma \sqrt{n}}.$$

How will this (deterministic, not random!) sequence of deviations behave, when compared to the probabilistic model?

Of course, for deterministic sequences, we cannot obtain cumulative or density distribution functions. Instead, we can only compute their histograms.

For example, to obtain a *normalized cumulative histogram* of the values $d_{10}(C, n)$ for $n \leq 10^9$, we let the variable $x$ run from $-3$ to $+3$ in steps of 0.1 (or 0.025), count all the indices $n$ such that $n \leq 10^9 \wedge d_{10}(C, n) \leq x$, and divide this count by $10^9$ (or $10^{10}$).

To obtain a *normalized density histogram* of the values $d_{10}(C, n)$ for $n \leq 10^9$, we let the variable $x$ run from $-2.9$ to $+3$ by step 0.1 (or 0.025), count all the indices $n$ such that $n \leq 10^9 \wedge x - step < d_{10}(C, n) \leq x$, and divide this count by $10^9$ (or $10^{10}$).

Normalization is necessary to make our histograms comparable to the corresponding cumulative and density distribution functions of the probabilistic model.

In Fig. 1, the density of standard normal distribution:

$$\phi(x) = \frac{d\,\Phi(x)}{dx} = \frac{1}{\sqrt{2\pi}} e^{-\frac{x^2}{2}}$$

is represented as a solid red curve. The blue bars represent then normalized density histogram of the values $d_{10}(\sqrt{2}, n)$ for $n \leq 10^9$.

In Fig. 1a, the cumulative distribution function $\Phi(x)$ of the standard normal distribution is represented as a smooth red line. The blue line represents the *normalized cumulative histogram* of the values $d_{10}(\sqrt{2}, n)$ for $n \leq 10^9$.

Naturally, both histograms convey the same information: the cumulative histogram can be obtained by summing the heights of the bars in the density histogram, while the density histogram can be derived by taking the differences between successive values in the cumulative histogram. Likewise, the red lines are related as the integral and derivative of each other.

As $n \to \infty$, how "should" these histograms behave in relation to the curves predicted by CLT in the probabilistic model? Should they tend to approach the theoretical curves?

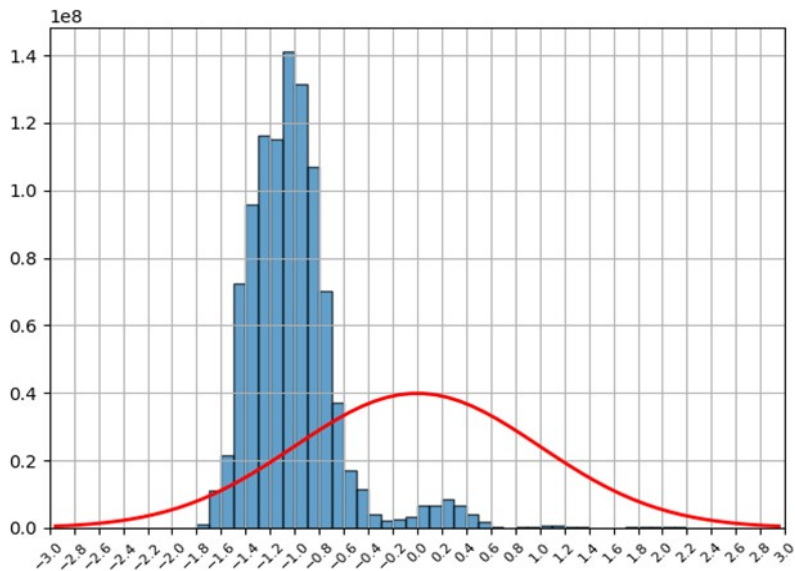

**Fig. 1**. The number $\sqrt{2}$: normalized density histogram of $d_{10}(\sqrt{2}, n)$, $n \leq 10^9$.
The red curve represents the density of standard normal distribution.

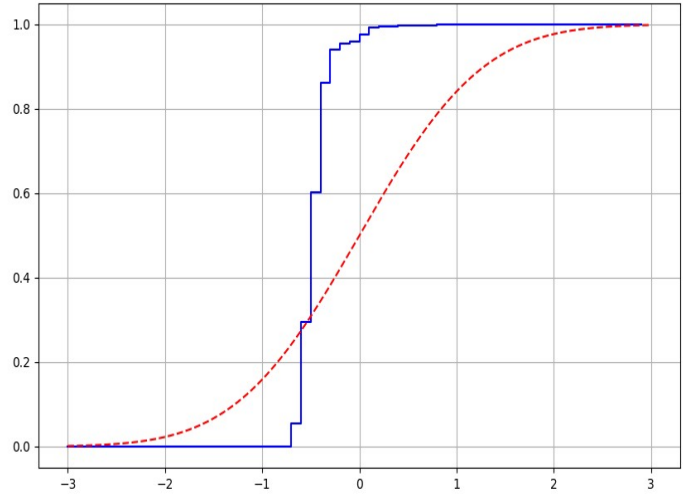

**Fig. 1a.** The number $\sqrt{2}$: normalized cumulative histogram of $d_{10}(\sqrt{2}, n)$, $n \le 10^9$.
The red line represents cumulative distribution function of standard normal distribution.

However, as seen in Fig. 1 and Fig. 1a, this is not the case – at least, not when considering the first billion decimal digits of the number $\sqrt{2}$.

## 4. Law of the Iterated Logarithm (LIL)

Belshaw and Borwein (2013) proposed applying the Law of the Iterated Logarithm (LIL) as a stronger test of the randomness of base-q representations of numbers.

In the probabilistic model described above, the LIL can be seen as a refinement of the CLT.

Let us once again consider $S_n^{(q)}$, the sum of the first $n$ base-$q$ random digits, and its centered and normalized deviation from its expected value $\mu n$:

$$d_n^{(q)} = \frac{S_n^{(q)} - \mu n}{\sigma \sqrt{n}}.$$

Now, let us refine this deviation by dividing it by $\sqrt{2 \log \log n}$ (where *log* denotes the natural logarithm),

$$\delta_n^{(q)} = \frac{d_n^{(q)}}{\sqrt{2 \log \log n}} = \frac{S_n^{(q)} - \mu n}{\sqrt{2 \sigma^2 n \log \log n}}.$$

*Note.* Recall that $\mu = \frac{q-1}{2}$; $\sigma^2 = \frac{q^2-1}{12}$.

According to LIL, in the probabilistic model, with probability 1, as $n \rightarrow \infty$, the values of $\delta_n^{(q)}$ oscillate infinitely between –1 and +1. More precisely, with probability 1:

$$\liminf_{n \to \infty} \delta_n^{(q)} = -1;$$

$$\limsup_{n \to \infty} \delta_n^{(q)} = +1.$$

Now, let us return to the digits of mathematical constants, and once again consider the base-$q$ representation of some mathematical constant $C$.

Let us consider the sequence of numbers derived in a "LIL-manner":

$$\delta_q(C, n) = \frac{d_q(C, n)}{\sqrt{2 \log \log n}} = \frac{S_q(C, n) - \mu n}{\sqrt{2 \sigma^2 n \log \log n}}.$$

How "should" this (deterministic!) sequence behave? Should it exhibit behavior as "predicted by LIL". Namely, should we expect that:

$$\liminf_{n \to \infty} \delta_q(C, n) = -1;$$

$$\limsup_{n \to \infty} \delta_q(C, n) = +1?$$

In other words, should $\delta_q(C, n)$, as a function of $n$, oscillate infinitely across the entire interval $(-1, +1)$. with the graph of $\delta_q(C, n)$ resembling a real oscillogram?

*Note*. In the probabilistic model, LIL holds with probability 1. However, in our case, we are dealing with a *deterministic* sequence – the digits of a mathematical constant $C$. Thus, $\delta_q(C, n)$ represents a completely fixed function of $n$.

Consequently, statements such as $\liminf_{n \to \infty} \delta_q(C, n) = -1$, or $\limsup_{n \to \infty} \delta_q(C, n) = +1$ are either *definitively true* or *false*, entirely independent of any probabilistic reasoning. (Nevertheless, in the worst case, such statements may turn out to be undecidable in ZFC.)

We applied this straightforward test to the first billion decimal digits of twelve mathematical constants: $\pi$, $e$, Euler–Mascheroni constant (γ), Apéry's constant $\zeta(3)$, quadratic irrationalities – square roots of 2 and 3 and the golden ratio, natural logarithms – *log 2, log 3, log 10*, Catalan's constant (G), and lemniscate constant ($\varpi$). Ten billion digits were examined in the case of $\pi$.

Unfortunately, actual calculations present a picture that differs significantly from the one predicted by LIL, at least for the first billion decimal digits (for ten billion of digits of $\pi$) – see the diagrams in the Appendix.

In particular, in the case of √2, as shown in the diagrams in the Appendix, the curve of $\delta_{10}(\sqrt{2}, n)$ oscillates within the interval $(-0.60, -0.15)$, starting from $n = 0.1 \cdot 10^9$ and continuing up to $n = 10^9$.

As a result, more than 90% of the values of

$$d_{10}(\sqrt{2},n)=\delta_{10}(\sqrt{2},n)\sqrt{2\log\log n}$$

fall between $-0.60\sqrt{2\log\log 10^9}\approx-1.48$ and $-0.15\sqrt{2\log\log 10^8}\approx-0.36$.

This conclusion is also reflected in the density histogram of $d_{10}(\sqrt{2},n)$ (Fig. 1), which is shifted to the left of the center 0, with most of the bar volume falling between $-1.48$ and $-0.36$.

Of course, it would be incorrect to assume that the behavior of the curve of $\delta_{10}(\sqrt{2},n)$ "explains" the shifted nature of the histograms of $d_{10}(\sqrt{2},n)$. In fact, given the limited number of digits considered, the curve of $\delta_{10}(\sqrt{2},n)$ conveys essentially the same information as the histograms of $d_{10}(\sqrt{2},n)$ shown in Fig. 1 and Fig. 1a.

*Note*. Incidentally, $\sqrt{2\log\log n}$ grows extremely slowly: for $n=10^8;10^9;10^{10}$ it takes on the values $2.41;2.46;2.50$, respectively.

## 5. The case of the number $\pi$

A situation similar to that of $\sqrt{2}$, was observed for six other constants: $\pi,\gamma,\zeta(3),\sqrt{3},\log 3,\bar{\omega}$, specifically, for those whose respective "LIL curves" of $\delta_{10}(C,n)$ tend to remain consistently above or below the zero level. Equivalently, the corresponding histograms of $d_{10}(C,n)$ are shifted away from the center 0.

Surprising, isn't it? One would expect the density histograms to become symmetric over time – roughly balanced on both sides of the center at 0. Does this not suggest that, for these seven constants (including $\sqrt{2}$), the future behavior of the curve of $\delta_{10}(C,n)$ will eventually shift – approaching the oscillations between –1 and +1, as predicted by LIL?

However, there is no indication of such a transition within the respective first billions of digits. In fact, at least in the case of $\pi$, this behavior remains unchanged even when the analysis is extended to ten billion digits.

What do the diagrams in Fig. 2, 2a, 2b, 2c reveal about the behavior of the corresponding curve of $\delta_{10}(\pi,n)$? In Fig. 2b we observe that for $n\leq 10^{10}$ most values of $d_{10}(\pi,n)$ fall between 0 and +1.6. Consequently, the majority of the corresponding values of $\delta_{10}(\pi,n)$ must fall within the interval from 0 and $\dfrac{1.6}{\sqrt{2\log\log 10^9}}\approx 0.65$. And this is confirmed by the behavior of the curve $\delta_{10}(\pi,n)$ from $n=10^9$ up to $n=10^{10}$.

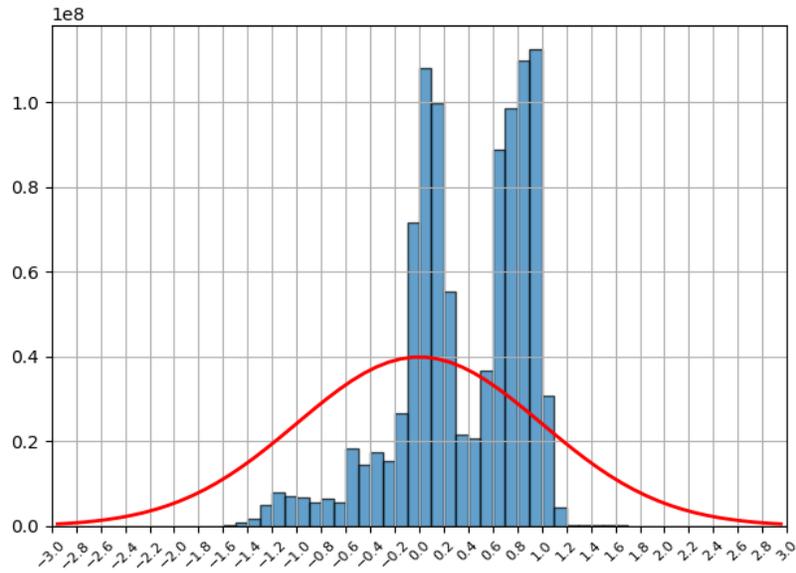

**Fig. 2.** The number $\pi$: normalized density histogram of $d_{10}(\pi, n)$, $n \leq 10^9$.
The red curve represents the density of standard normal distribution.

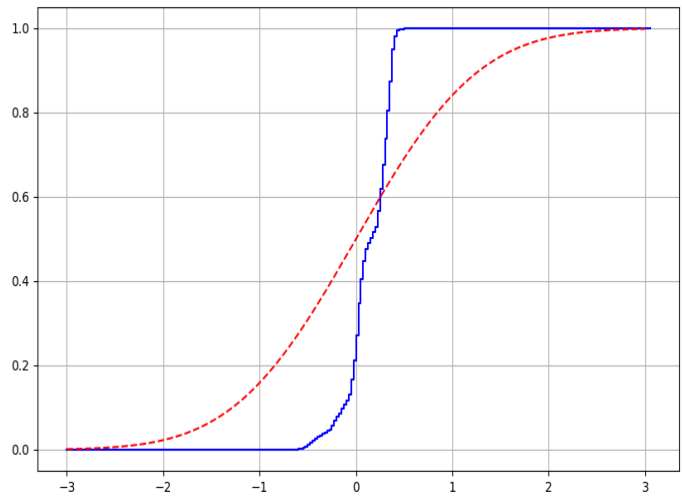

**Fig. 2a.** The number $\pi$: normalized cumulative histogram of $d_{10}(\pi, n)$, $n \leq 10^9$ (bin 0.025)
The red line represents cumulative distribution function of standard normal distribution.

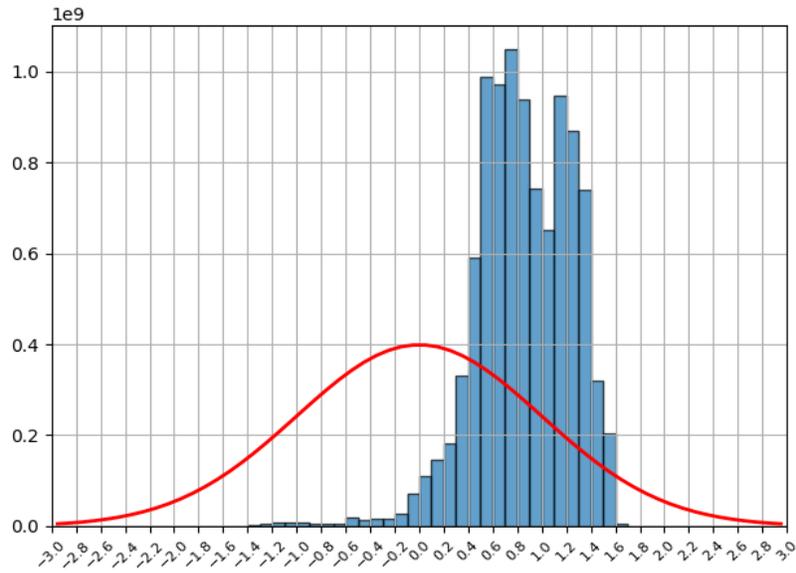

**Fig. 2b.** The number $\pi$: normalized density histogram of $d_{10}(\pi, n)$, $\boldsymbol{n \leq 10^{10}}$.
The red curve represents the density of standard normal distribution.

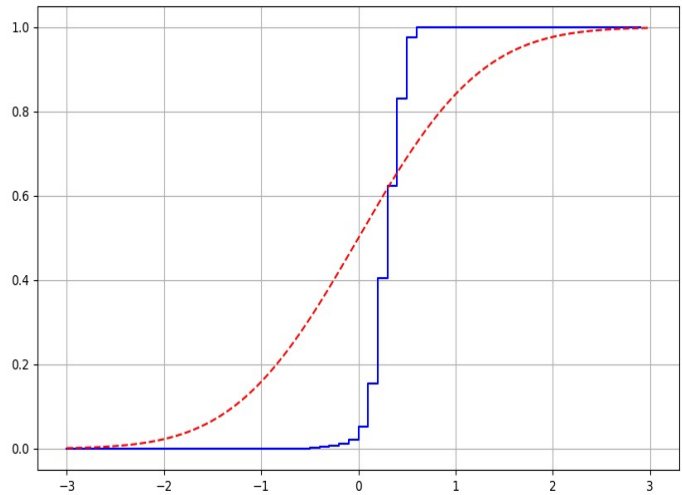

**Fig. 2c.** The number $\pi$: normalized cumulative histogram of $d_{10}(\pi, n)$, $\boldsymbol{n \leq 10^{10}}$ (bin 0.1)
The red line represents cumulative distribution function of standard normal distribution.

## 6. The case of the number *e*

The first version of this paper (Podnieks 2014) contains a conjecture: $\lim_{n \to \infty} \delta_{10}(C, n) = 0$, seemingly confirmed by the first $10^7$ decimal digits of the constants $\pi, e, \sqrt{2}$. However, as shown in the diagrams presented in the Appendix, this stronger conjecture fails shortly after $10^7$ digits – except (perhaps) for the digits of $e, \frac{1+\sqrt{5}}{2}, \log 2, \log 10$.

These constants appear to exhibit oscillations of $\delta_{10}(C, n)$ around 0, and the overall tendency does seem to be one of convergence:

$$\lim_{n \to \infty} inf \, \delta_{10}(C, n) = 0.$$

See also Fig. 3, 3a and 3b for the number e.

*Note.* Catalan's G also shows oscillations around 0, but without any clear tendency toward convergence.

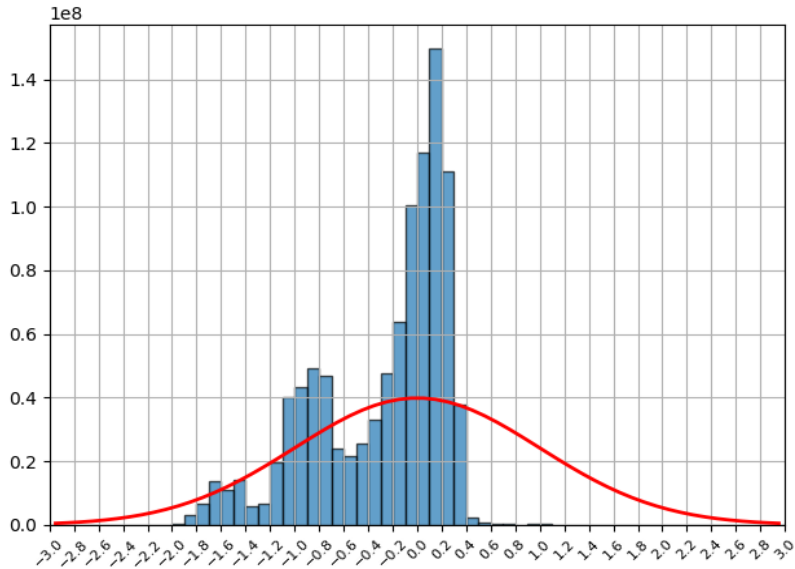

**Fig. 3.** The number *e*: normalized density histogram of $d_{10}(e, n), n \leq 10^9$.
The red curve represents the density of standard normal distribution.

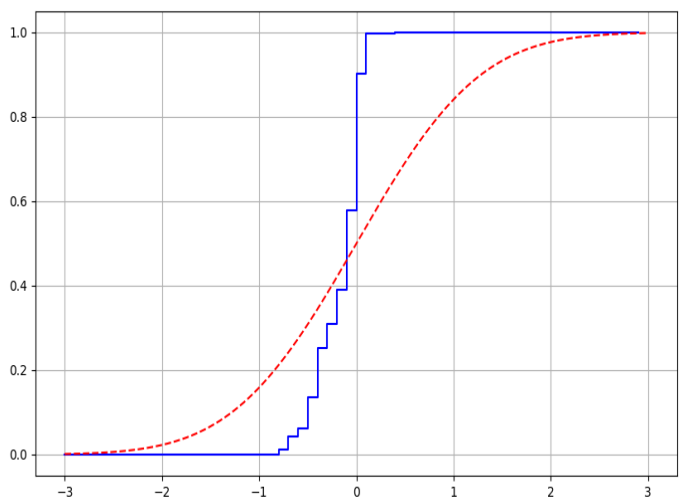

**Fig. 3a.** The number $e$: normalized cumulative histogram of $d_{10}(e, n)$, $n \le 10^9$.
The red line represents cumulative distribution function of standard normal distribution.

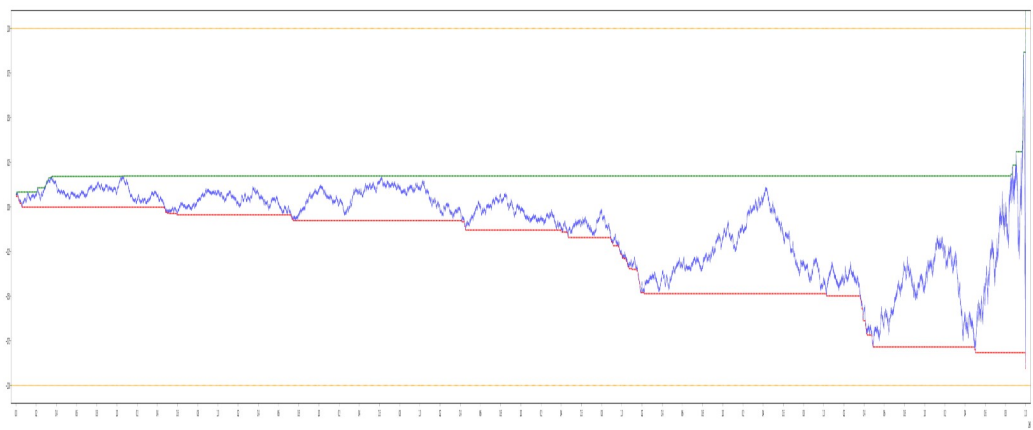

**Fig. 3b.** The number $e$: simulated *lim inf* and *lim sup* of $\delta_{10}(e, n)$.
From left to right, $n$ is running from $10^9$ down to 1.

## 7. Conclusion

According to a popular belief, the decimal digits of mathematical constants such as π behave like statistically independent random variables, each taking the values 0, 1, 2, 3, 4, 5, 6, 7, 8, and 9 with equal probability of 1/10. If this is the case, then, in particular, the decimal representations of these constants should tend to satisfy the *central limit theorem* (CLT) and the *law of the iterated logarithm* (LIL).

The paper presents the results of a direct statistical analysis of the decimal representations of 12 mathematical constants with respect to the central limit theorem (CLT) and the law of the iterated logarithm (LIL). The analyzed constants include π, e, the Euler–Mascheroni constant (γ), Apéry's constant (ζ(3)), quadratic irrationals such as the square roots of 2 and 3, and the golden ratio. Additionally, natural logarithms – *log 2*, *log 3*, and *log 10*, Catalan's constant (G), and the lemniscate constant ($\overline{\omega}$) were examined.

The first billion digits of each constant were analyzed, with ten billion digits examined in the case of π. Within these limits, no evidence was found to suggest that the digits of these constants satisfy CLT or LIL.

So, do the first billion digits of π, or of another "naturally occurring" mathematical constant, represent a "truly random" body of data?

Consider, for example, the first billion of digits of $\sqrt{2}$ or $e$ , and look at Fig. 1 and Fig. 3. Both, $d_{10}(\sqrt{2},n)$ and $d_{10}(e,n)$ are less than 0.6 "almost surely". Thus, for almost all $n \leq 10^9$, the sum of digits $S_{10}(e,n)$ is less than $4.5n + 0.6\sqrt{8.25n}$, which contradicts the 68–95–99.7% rule (a consequence of CLT), according to which approximately $(100-68)/2 = 16\%$ of these sums should exceed $4.5n + \sqrt{8.25n}$.

In a similar way, we can conclude that none of the examined constants generate their first billion digits in a "truly random" way.

Therefore, all these billions of digits cannot be regarded as "randomly" generated bodies of data.

If we were to extend the calculations well beyond $10^9$ digits, would we observe oscillations of the curve $\delta_{10}(C,n)$ between –1 and +1? And would the histograms of $d_{10}(C,n)$ become symmetric?

So far, we have done this only for the number $\pi$ – up to its first ten billion digits. We did not observe any signs that the $\delta_{10}(\pi,n)$ curve might change its behavior: starting from $n = 1.9 \cdot 10^9$, it appears to oscillate within the interval $(+0.24, +0.60)$, rather than within the "expected" (–1, +1) range.

To conclude: The first billion decimal digits of the "naturally occurring" mathematical constants considered do not provide a truly random body of data. Beyond this "curtain" of $10^9$ digits, the situation appears unchanged – at least up to the tenth billion digit of $\pi$.

Consequently, it does not seem worthwhile to invest computational resources in exploring $10^{10}$ or even $10^{10}$ decimal digits of such constants. Does exploring $10^{12}$ digits look any more promising?

## Appendix

The first billion decimal digits of each constant examined (except $\pi$) were downloaded from the Internet Archive (2016). The first ten billion decimal digits of $\pi$ were obtained from a website maintained by Jack Giffin.

To verify the accuracy of these sources, portions of the digits of $\pi$ and $\sqrt{2}$ were recalculated using the Bailey–Borwein–Plouffe (BBP) algorithm (Bailey, 2006, for $\pi$) and the Newton–Raphson method (for $\sqrt{2}$), respectively.

The calculation of the values of $\delta_{10}(C,n)$ using the formula:

$$\delta_{10}(C,n) = \frac{d_{10}(C,n)}{\sqrt{2\log\log n}} = \frac{S_{10}(C,n) - 4.5\,n}{\sqrt{2 \cdot 8.25\,n\log\log n}}$$

did not pose any precision issues, as the sums $S_{10}(C,n)$ were accumulated as integers.

The diagrams below are organized as follows: each diagram displays the values of $\delta_{10}(C,n)$ for a consecutive block of 100 million values of $n$.

For the Python and C++ source codes, as well as the diagrams in separate PNG files, see Roba (2024a).

For the full set of computed $\delta_{10}(C,n)$ values, refer to Roba (2024b).

The number $\pi$ (ten billions of digits)

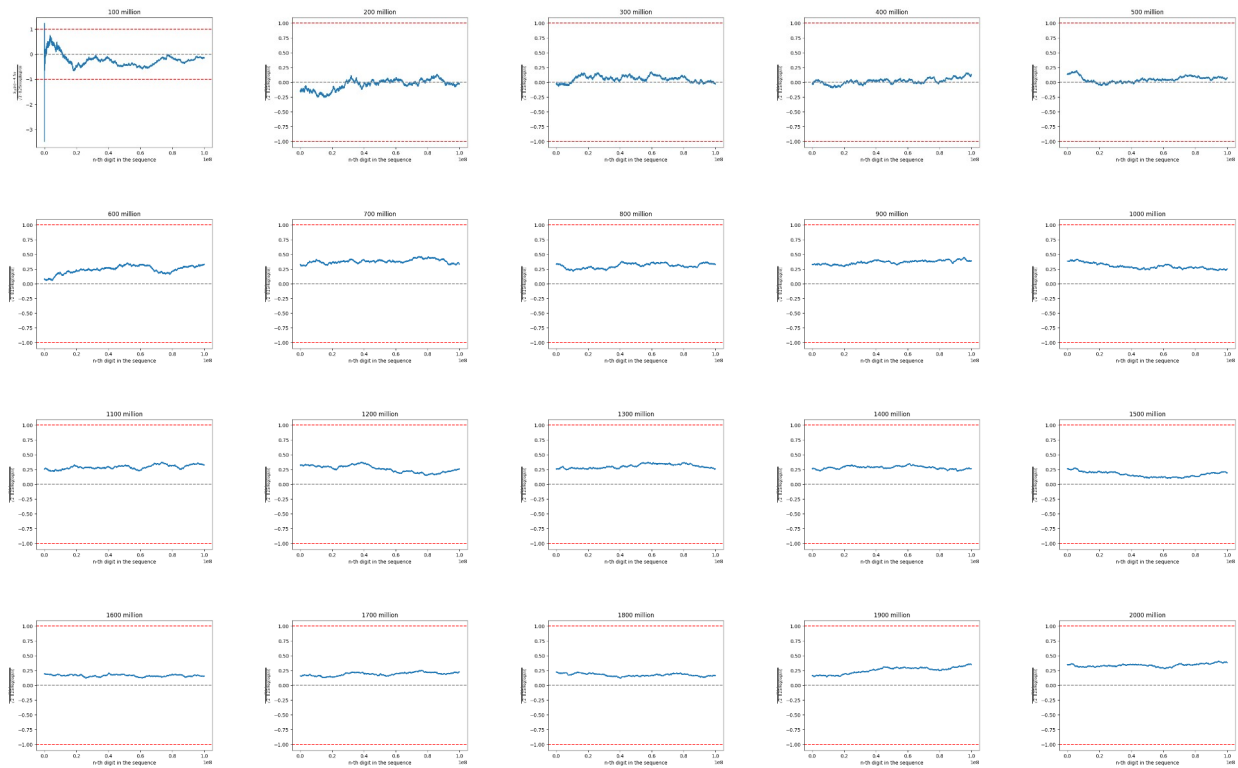

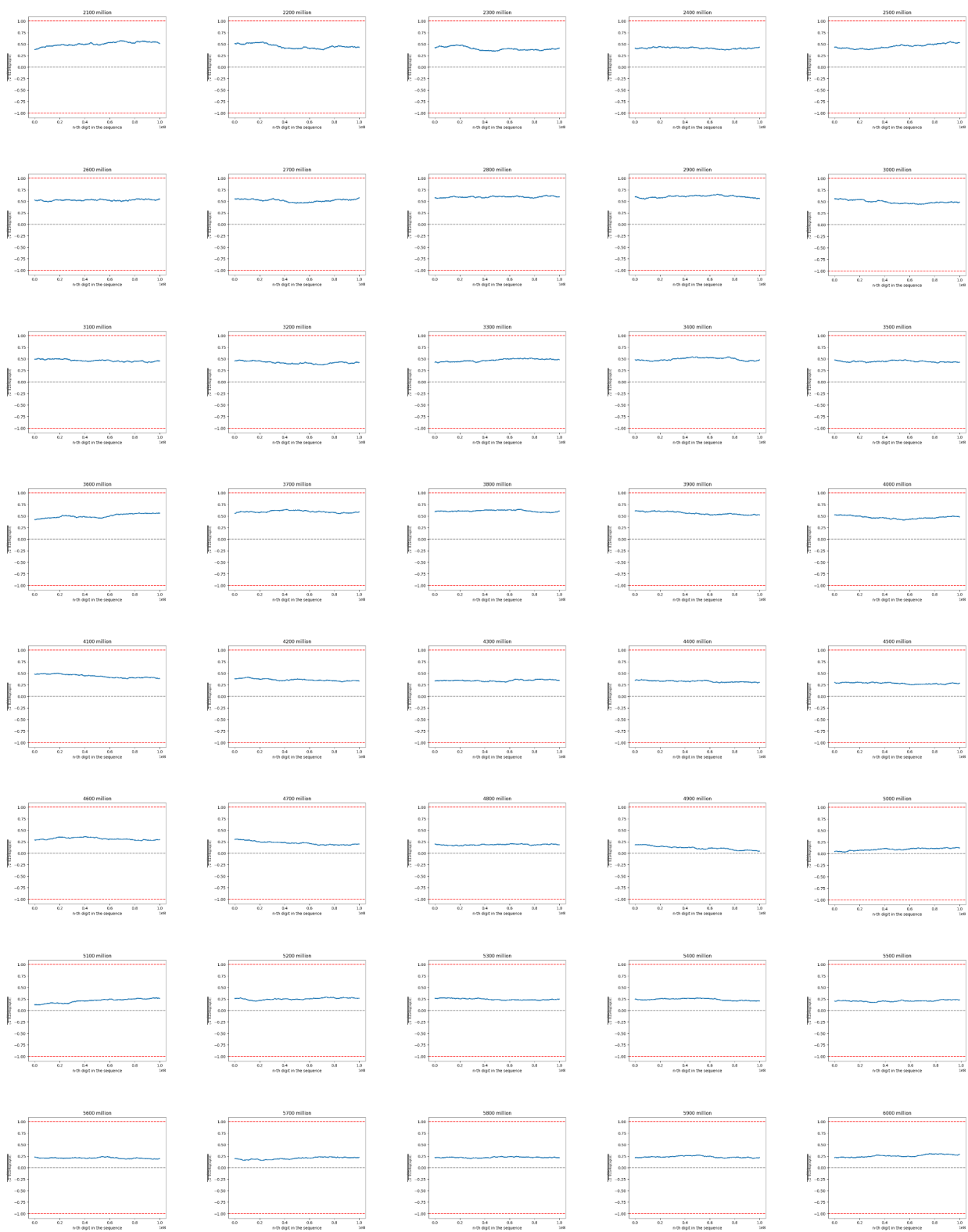

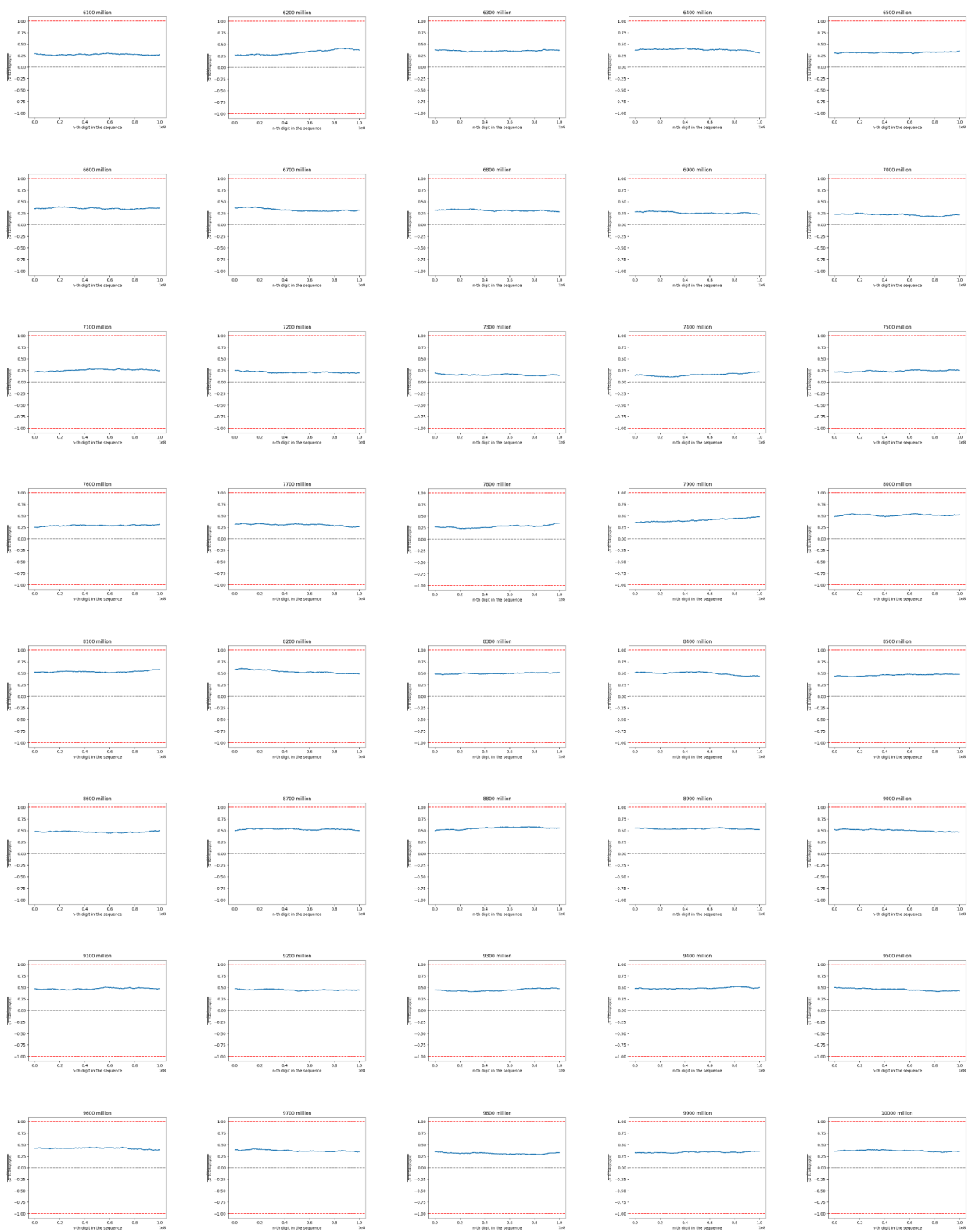

Euler number e (one billion of digits)

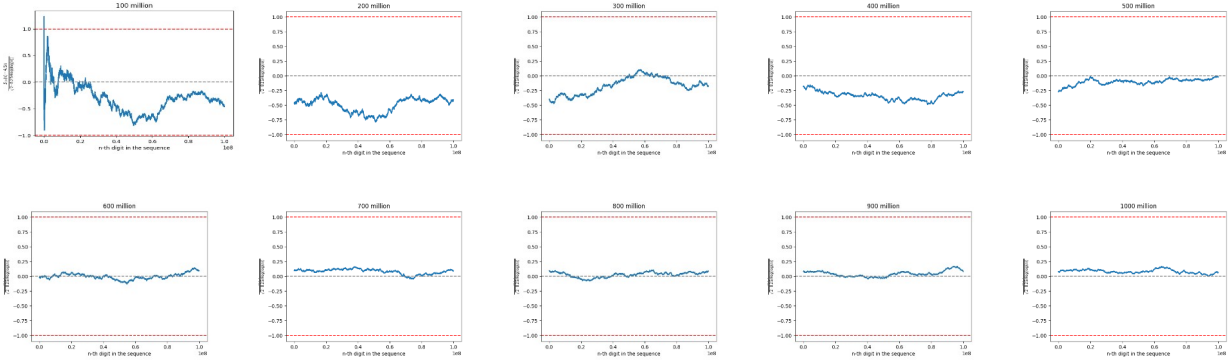

Euler-Mascheroni constant $\gamma$ (one billion of digits)

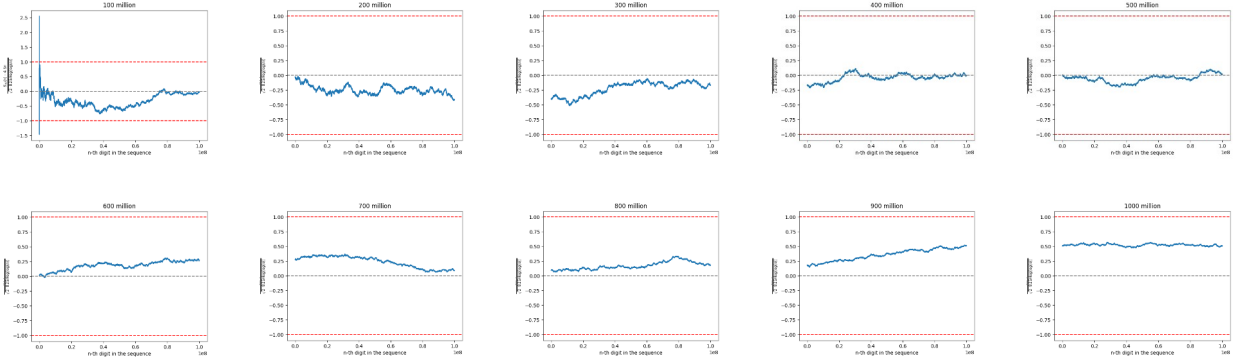

Apéry's constant $\zeta(3)$ (one billion of digits)

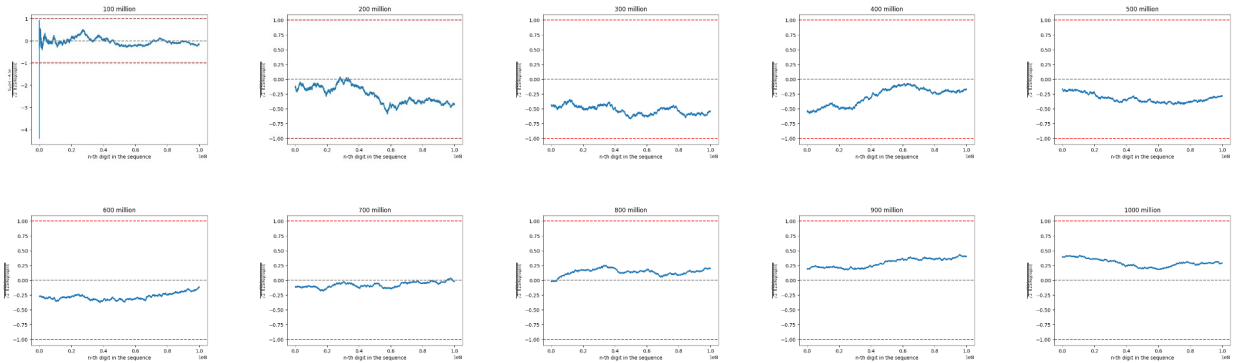

$\sqrt{2}$ (one billion of digits)

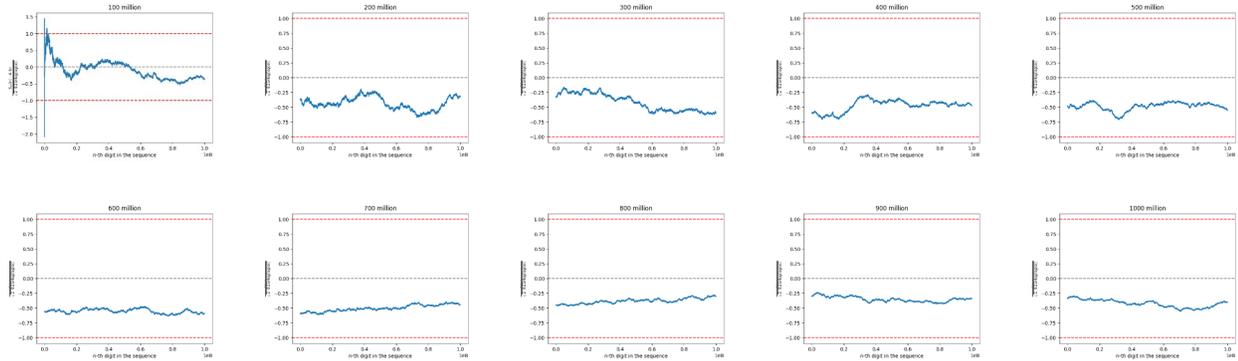

$\sqrt{3}$ (one billion of digits)

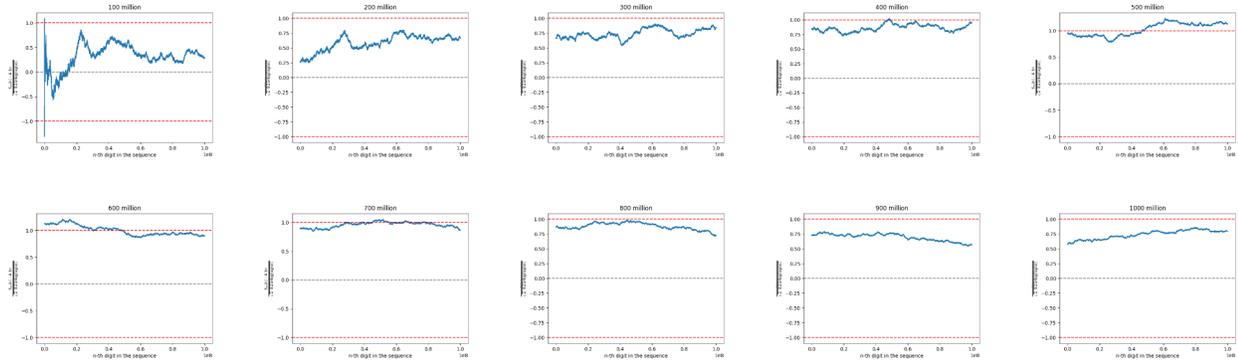

Golden ratio $\frac{1+\sqrt{5}}{2}$ (one billion of digits)

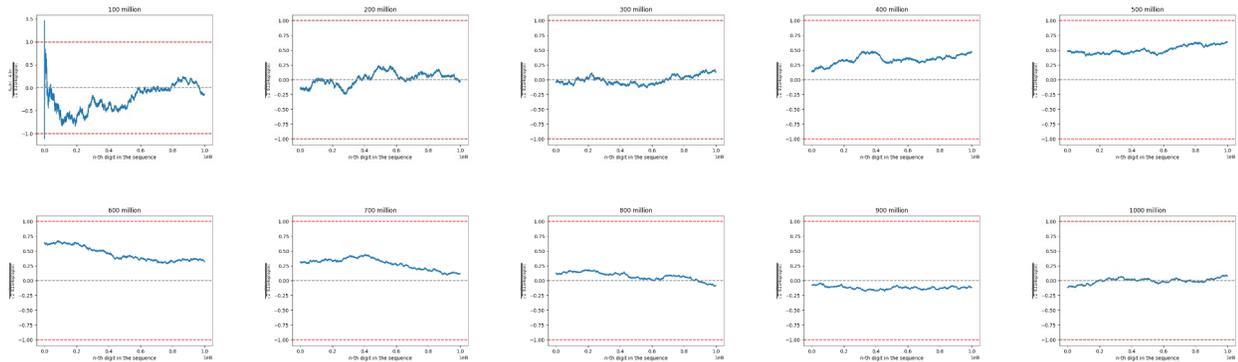

log 2 (one billion of digits)

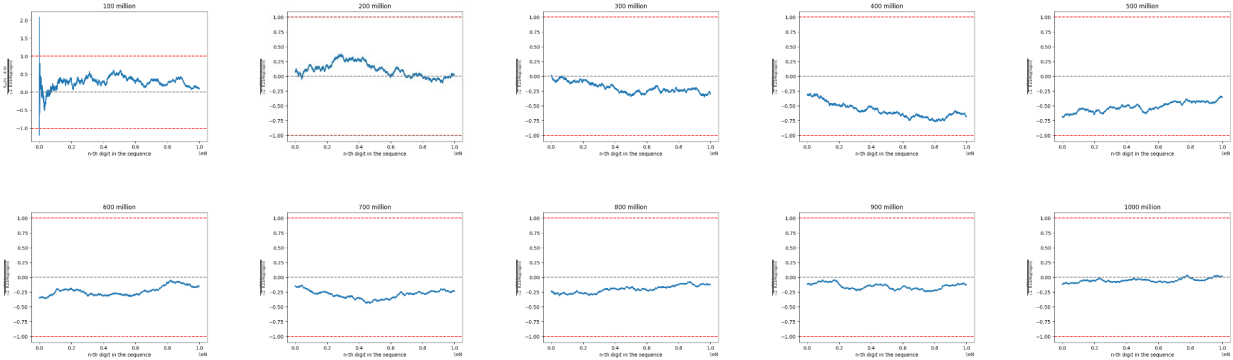

log 3 (one billion of digits)

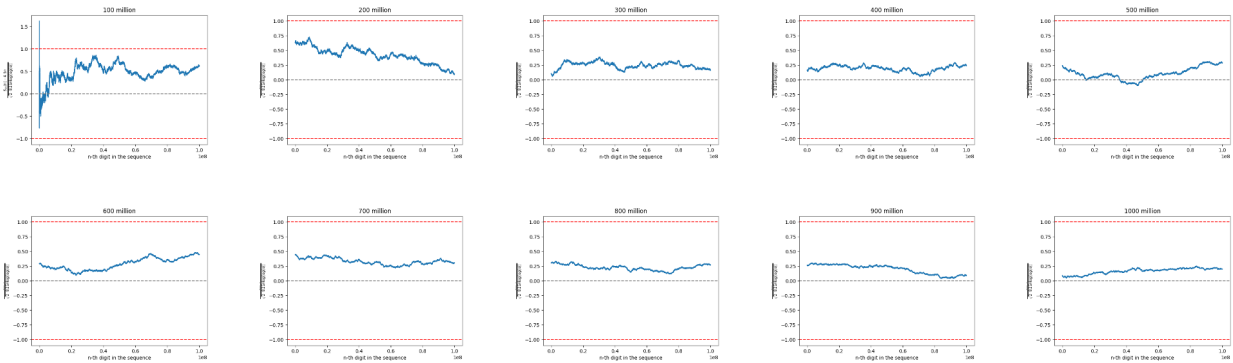

log 10 (one billion of digits)

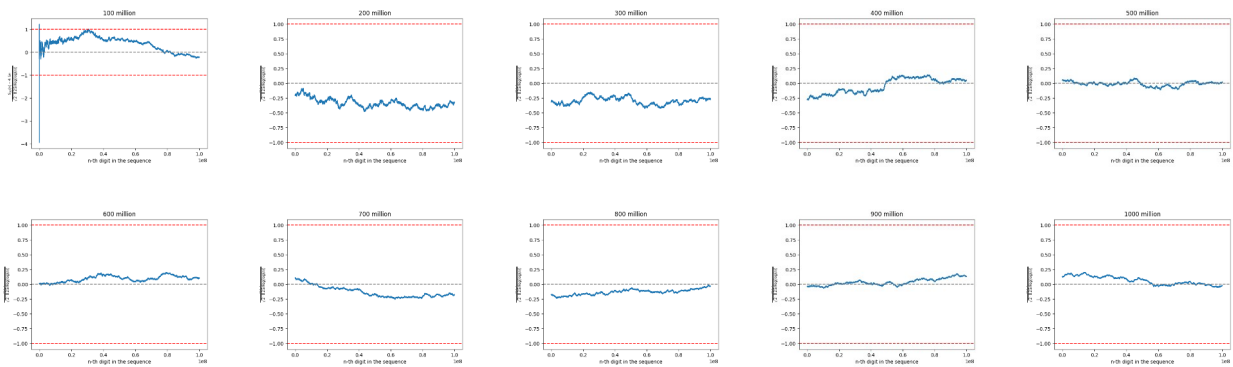

Catalan's constant $G$ (one billion of digits)

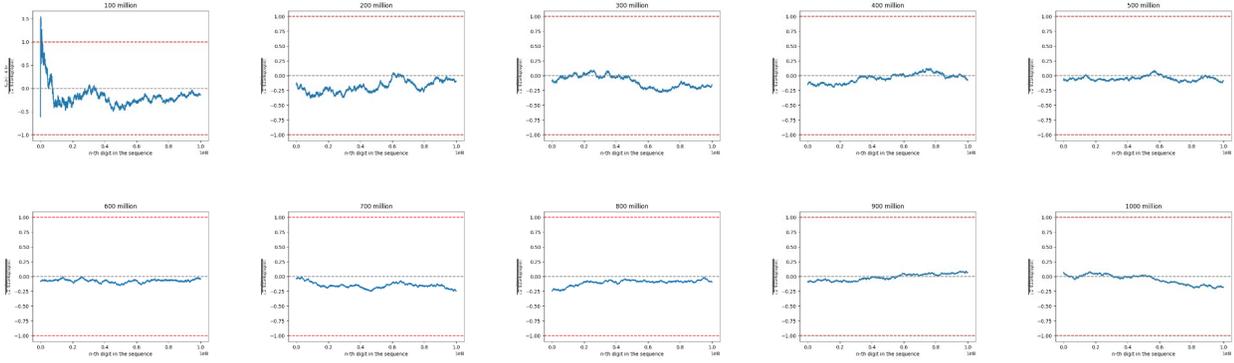

Lemniscate constant $\bar{\omega}$ (one billion of digits)

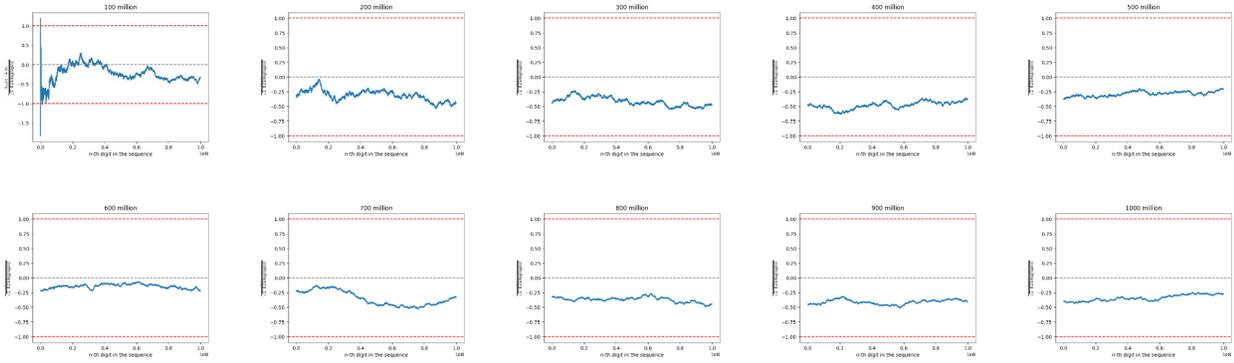